\documentclass[a4paper,12pt]{amsart}
\newtheorem{thm}{Theorem}[section]
\newtheorem{defi}[thm]{Definition}
\newtheorem{conj}[thm]{Conjecture}

\usepackage{amssymb}
\def\BA{{\mathbf{A}}}
\def\BC{{\mathbf{C}}}
\def\uc{\operatorname{\underline{c}}\nolimits}
\def\eps{\varepsilon}
\def\Gg{{\mathfrak{g}}}
\def\Gal{\operatorname{Gal}\nolimits}
\def\BH{{\mathbf{H}}}
\def\ie{{\em i.e.}}
\def\Ind{\operatorname{Ind}\nolimits}
\def\Irr{\operatorname{Irr}\nolimits}
\def\iso{\buildrel \sim\over\to}
\def\GL{\operatorname{GL}\nolimits}
\def\CO{{\mathcal{O}}}
\def\Gp{{\mathfrak{p}}}
\def\CP{{\mathcal{P}}}
\def\CQ{{\mathcal{Q}}}
\def\Gq{{\mathfrak{q}}}
\def\BR{{\mathbf{R}}}
\def\CR{{\mathcal{R}}}
\def\Gr{{\mathfrak{r}}}
\def\CS{{\mathcal{S}}}
\def\GS{{\mathfrak{S}}}
\def\Spec{\operatorname{Spec}\nolimits}
\def\BZ{{\mathbf{Z}}}

\title{Calogero-Moser versus Kazhdan-Lusztig cells}
\author{C\'edric Bonnaf\'e and Rapha\"el Rouquier}
\address{C\'edric Bonnaf\'e~: Universit\'e Montpellier 2,
Institut de Math\'ematiques et de Mod\'elisation de Montpellier,
Case Courrier 051, Place Eug\`ene Bataillon, 34095 Montpellier Cedex,
FRANCE}
\email{cedric.bonnafe@math.univ-montp2.fr}
\address{Rapha\"el Rouquier~: Mathematical Institute,
University of Oxford, 24-29 St Giles', Oxford, OX1 3LB, UK
and Department of Mathematics, UCLA, Box 951555,
Los Angeles, CA 90095-1555, USA}
\email{rouquier@maths.ox.ac.uk}

\date\today
\begin{document}
\maketitle
\section{Introduction}
In \cite{KaLu}, Kazhdan and Lusztig developed a combinatorial theory
associated with Coxeter groups. They defined in particular partitions
of the group in left and two-sided cells.
For Weyl groups, these have a representation theoretic interpretation
in terms of primitive ideals, and they
play a key role in Lusztig's description of unipotent
characters for finite groups of Lie type \cite{Lu3}.
Lusztig generalized this theory to Hecke algebras of Coxeter groups with
unequal parameters \cite{Lu2,Lu4}.

We propose a definition of left cells and two-sided cells for
complex reflection groups, based on ramification theory for Calogero-Moser
spaces. These spaces have been defined via rational Cherednik algebras by
Etingof and Ginzburg \cite{EtGi}. We conjecture
that these coincide with Kazhdan-Lusztig cells, for real reflection groups.
Counterparts of families of irreducible characters have been studied by
Gordon and Martino \cite{GoMa}, and we provide here a version of left
cell representations.
The Calogero-Moser cells are studied in detail in \cite{BoRou}.

\section{Calogero-Moser spaces and cells}
\subsection{Rational Cherednik algebras at $t=0$}
Let us recall some constructions and results from \cite{EtGi}.
Let $V$ be a finite-dimensional complex vector space and
$W$ a finite subgroup of $\GL(V)$. Let $\CS$ be the set of reflections of
$W$, \ie, elements $g$ such that $\ker(g-1)$ is a hyperplane.
We assume that $W$ is a reflection group, \ie, it is generated by $\CS$.

We denote by $\CS/\sim$ the quotient of $\CS$ by the conjugacy action of $W$ and
we let $\{\uc_s\}_{s\in \CS/\sim}$ be a set of indeterminates.
We put $A=\BC[\BC^{\CS/\sim}]=\BC[\{\uc_s\}_{s\in \CS/\sim}]$.
Given $s\in \CS$, let $v_s\in V$ (resp. $\alpha_s\in V^*$) be an eigenvector
for $s$ associated to the non-trivial eigenvalue.

The $0$-rational Cherednik algebra $\BH$ is the quotient of
$A\otimes T(V\oplus V^*)\rtimes W$ by the relations
$$[x,x']=[\xi,\xi']=0$$
$$[\xi,x]=\sum_{s\in \CS} \uc_s \frac{\langle x,\alpha_s\rangle\cdot
\langle v_s,\xi\rangle}{\langle v_s,\alpha_s\rangle}s
\ \text{ for }x,x'\in V^* \text{ and }\xi,\xi'\in V.$$

We put $Q=Z(\BH)$ and $P=A\otimes S(V^*)^W\otimes S(V)^W\subset Q$. The ring
$Q$ is normal.
It is a free $P$-module of rank $|W|$.

\subsection{Galois closure}
Let $K=\mathrm{Frac}(P)$ and $L=\mathrm{Frac}(Q)$. Let $M$ be a Galois closure
of the extension $L/K$ and $R$ the integral closure of $Q$ in $M$.
Let $G=\Gal(M/K)$ and $H=\Gal(M/L)$. Let 
$\CP=\Spec P=\BA_\BC^{\CS/\sim}\times V/W\times V^*/W$, $\CQ=\Spec Q$ the
Calogero-Moser space, and $\CR=\Spec R$.

We denote by $\pi:\CR\to\CQ$ the quotient by $H$, and by
$\Upsilon:\CQ\to\CP$ and
$\phi:\CP\to \BA_\BC^{\CS/\sim}$ the canonical maps. We put
$p=\Upsilon\pi:\CR\to\CP$ the quotient by $G$.

\subsection{Ramification}
Let $\Gr\in\CR$ be a prime ideal of $R$. We denote by $D(\Gr)\subset G$ its
decomposition group and by $I(\Gr)\subset D(\Gr)$ its inertia group.

We have a decomposition into irreducible components
$$\CR\times_\CP \CQ=\bigcup_{g\in G/H}
\CO_g, \text{ where }
\CO_g=\{(x,\pi(g^{-1}(x)))| x\in \CR\}$$
inducing a decomposition into irreducible components
$$V(\Gr)\times_\CP \CQ=\coprod_{g\in I(\Gr)\setminus G/H}\CO_g(\Gr),
\text{ where }
\CO_g(\Gr)=\{(x,\pi(g^{-1}g'(x)))| x\in V(\Gr),\ g'\in I(\Gr)\}.$$

\subsection{Undeformed case}
Let $\Gp_0=\phi^{-1}(0)=\sum_{s\in \CS/\sim}P \uc_s$.
We have
$P/\Gp_0=\BC[V\oplus V^*]^{W\times W}$,
$Q/\Gp_0 Q=\BC[V\oplus V^*]^{\Delta W}$, where $\Delta(W)=\{(w,w)| w\in W\}\subset
W\times W$. A Galois closure of the extension of
$\BC(\Gp_0 Q)=\BC(V\oplus V^*)^{\Delta W}$ over 
$\BC(\Gp_0)=\BC(V\oplus V^*)^{W\times W}$ is
$\BC(V\oplus V^*)^{\Delta Z(W)}$.

\smallskip
Let $\Gr_0\in\CR$ above $\Gp_0$.
Since $\Gp_0 Q$ is prime, we have
$G=D(\Gr_0)H=HD(\Gr_0)$ and $I(\Gr_0)=1$.
Fix an isomorphism
$\iota:\BC(\Gr_0)\iso \BC(V\oplus V^*)^{\Delta Z(W)}$ extending the canonical
isomorphism of $\BC(\Gp_0 Q)$ with $\BC(V\oplus V^*)^{\Delta W}$.

The application $\iota$
induces an isomorphism $D(\Gr_0)\iso (W\times W)/\Delta Z(W)$,
that restricts to an isomorphism $D(\Gr_0)\cap H\iso \Delta W/\Delta Z(W)$.
This provides a bijection
$G/H\iso (W\times W)/\Delta W$. Composing with the inverse of the bijection
$W\iso (W\times W)/\Delta W,\ w\mapsto (w,1)$, we obtain a bijection
$G/H\iso W$.

\smallskip
 From now on, we identify the sets $G/H$ and $W$ through this
bijection.
 Note that this bijection depends on the choices of $\Gr_0$ and of $\iota$.
 Since $M$ is the Galois closure of $L/K$, we have
$\bigcap_{g\in G}H^g=1$, hence the left action of $G$ on $W$ induces
an injection $G\subset\GS(W)$.

\subsection{Calogero-Moser cells}
\label{se:CMcells}
\begin{defi}
Let $\Gr\in\CR$. 
The $\Gr$-cells of $W$ are the orbits of $I(\Gr)$ in its action on $W$.
\end{defi}

Let $c\in\BA_\BC^{\CS/\sim}$. Choose $\Gr_c\in\CR$ with
$\overline{p(\Gr_c)}=\bar{c}\times 0\times 0$. The $\Gr_c$-cells are called
the
{\em two-sided Calogero-Moser $c$-cells} of $W$.
Choose now $\Gr_c^{\mathrm{left}}\in\CR$ contained in $\Gr_c$ with
$\overline{p(\Gr_c^{\mathrm{left}})}=\bar{c}\times V/W\times 0\in\CP$.
The $\Gr_c^{\mathrm{left}}$-cells are called the
{\em left Calogero-Moser $c$-cells} of $W$. We have $I(\Gr_c^{\mathrm{left}})
\subset I(\Gr_c)$. Consequently, every left cell is contained in a unique
two-sided cell.

\smallskip
The map sending $w\in W$ to $\pi(w^{-1}(\Gr_c))$
induces a bijection from the set of two-sided cells to 
$\Upsilon^{-1}(c\times 0\times 0)$.

\subsection{Families and cell multiplicities}

Let $E$ be an irreducible representation of $\BC[W]$. We extend it to a
representation of $S(V)\rtimes W$ by letting $V$ act by $0$. Let
$$\Delta(E)=e\cdot\Ind_{S(V)\rtimes W}^{\BH}(A\otimes_{\BC}E), 
\text{ where }
e=\frac{1}{|W|}\sum_{w\in W}w,$$
be the spherical Verma module associated with $E$. It is a
$Q$-module.

\smallskip
Let $c\in\BA_\BC^{\CS/\sim}$ and
let $\Delta^{\!\mathrm{left}}(E)=(R/{\Gr_c^{\mathrm{left}}})\otimes_P \Delta(E)$.

\begin{defi}
Given $\Gamma$ a left cell, we define the cell multiplicity 
$m_\Gamma(E)$ of $E$ as
the multiplicity
of $\Delta^{\!\mathrm{left}}(E)$ at the component $\CO_\Gamma(\Gr^{\mathrm{left}}_c)$.
\end{defi}

Note that
$\sum_{\Gamma} m_\Gamma(E)\cdot [\CO_\Gamma(\Gr^{\mathrm{left}}_c)]$ is the
support cycle of $\Delta^{\!\mathrm{left}}(E)$.

\smallskip
There is a unique two-sided cell $\Lambda$ containing all left cells
$\Gamma$ such that $m_{\Gamma}(E)\not=0$. Its image in $\CQ$ is the unique
$\Gq\in\Upsilon^{-1}(c\times 0\times 0)$ such that 
$(Q/\Gq)\otimes_Q \Delta(E)\not=0$. The corresponding map
$\Irr(W)\to \Upsilon^{-1}(c\times 0\times 0)$ is surjective, and its
fibers are the {\em Calogero-Moser
families} of $\Irr(W)$, as defined by Gordon \cite{Go1}.

\subsection{Dimension $1$}
Let $V$ be a one-dimensional complex vector space, let $d\ge 2$ and let
$W$ be the group of $d$-th roots of unity acting on $V$.
Let $\zeta=\exp(2i\pi/d)$, let $s=\zeta\in W$ and
$\underline{c}_i=\underline{c}_{s^i}$ for $1\le i\le d-1$. We
have $A=\BC[\underline{c}_1,\ldots,\underline{c}_{d-1}]$ and
$$\BH=A\langle x,\xi,s| sxs^{-1}=\zeta^{-1} x,\ s\xi s^{-1}=\zeta\xi
\text{ and }[\xi,x]=\sum_{i=1}^{d-1}\underline{c}_is^i\rangle.$$
Let $\mathrm{eu}=\xi x-\sum_{i=1}^{d-1}(1-\zeta^i)^{-1}\underline{c}_is^i$. We have
$P=A[x^d,\xi^d]$ and $Q=A[x^d,\xi^d,\mathrm{eu}]$.
Define $\underline{\kappa}_1,\ldots,\underline{\kappa}_d=\underline{\kappa}_0$
by $\underline{\kappa}_1+\cdots+\underline{\kappa}_d=0$ and
$\sum_{i=1}^{d-1}\underline{c}_is^i=\sum_{i=0}^{d-1}(\underline{\kappa}_i-
\underline{\kappa}_{i+1})\eps_i$, where $\eps_i=\frac{1}{d}\sum_{j=0}^{d-1}\zeta^{ij}s^j$.
We have $A=\BC[\underline{\kappa}_1,\ldots,\underline{\kappa}_d]/(\underline{\kappa}_1+\cdots+\underline{\kappa}_d)$.

\smallskip
The normalization of the Galois closure is described as follows.
There is an isomorphism of $A$-algebras
$$A[X,Y,Z]/\bigl(XY-\prod_{i=1}^d(Z-\underline{\kappa}_i)\bigr)\iso Q,\ X\mapsto x^d,\ Y\mapsto \xi^d
\text{ and }Z\mapsto \mathrm{eu}.$$
We have an isomorphism of $A$-algebras
$$A[X,Y,\lambda_1,\ldots,\lambda_d]/\bigl(e_1(\lambda)=
e_1(\underline{\kappa}),\ldots,
e_{d-1}(\lambda)=e_{d-1}(\underline{\kappa}),e_d(\lambda)=
e_d(\underline{\kappa})+(-1)^{d+1}XY\bigr)\iso R$$
where $Z=\lambda_d$ and where $e_i$ denotes the $i-th$ elementary symmetric
function. We have $G=\GS_d$, acting by permuting the $\lambda_i$'s, and
$H=\GS_{d-1}$. 

Let $\Gp_0=(\underline{\kappa}_1,\ldots,\underline{\kappa}_d)\in\Spec P$ and
$\Gr_0=(\underline{\kappa}_1,\ldots,\underline{\kappa}_d,
\lambda_1-\zeta \lambda_d,\ldots,\lambda_{d-1}-\zeta^{d-1}\lambda_d)
\in\Spec R$. We have $D(\Gr_0)=\langle (1,2,\ldots,d)\rangle\subset\GS_d$
and $\BC(\Gr_0)=\BC(X,Y,\lambda_d=\sqrt[d]{XY})=\BC(X,Y,Z=\sqrt[d]{XY})$.
The composite bijection $D(\Gr_0)\iso G/H\iso W$ is an isomorphism of groups
given by $(1,\ldots,d)\mapsto s$.

\smallskip
Fix $c\in\BC^{d-1}$ and let $\kappa_1,\ldots,\kappa_d\in\BC$ corresponding to
$c$. Consider $\Gr=\Gr_c$ or $\Gr_c^{\mathrm{left}}$ as in \S \ref{se:CMcells}. Then
$I(\Gr)$ is the subgroup of $\GS_d$ stabilizing $(\kappa_1,\ldots,
\kappa_d)$. The left $c$-cells coincide with the 
two-sided $c$-cells and two elements $s^i$ and $s^j$ are in the same cell
if and only if $\kappa_i=\kappa_j$. Finally, the multiplicity
$m_\Gamma(\det^j)$ is $1$ if $s^j\in\Gamma$ and $0$ otherwise.

\section{Coxeter groups}

\subsection{Kazhdan-Lusztig cells}
Following Kazhdan-Lusztig \cite{KaLu} and Lusztig \cite{Lu2,Lu4}, let us
recall the construction of cells.

We assume here $V$ is the complexification of a real vector space $V_\BR$
acted on by $W$. We choose a connected component $C$ of
$V_\BR-\bigcup_{s\in \CS}\ker(s-1)$ and we denote by $S$ the set of 
$s\in\CS$ such that $\ker(s-1)\cap \bar{C}$ has codimension $1$ in $\bar{C}$.
This makes $(W,S)$ into a Coxeter group, and we denote by $l$ the length
function.

\smallskip
Let $\Gamma$ be a totally ordered free abelian group and
let $L:W\to\Gamma$ be a weight function, \ie, a function
such that $L(ww')=L(w)+L(w')$ if $l(ww')=l(w)+l(w')$. We denote
by $v^\gamma$ the element of the group algebra $\BZ[\Gamma]$ corresponding
to $\gamma\in\Gamma$.

\smallskip
We denote by $H$ the Hecke algebra of $W$: this is the $\BZ[\Gamma]$-algebra
generated by elements $T_s$ with $s\in S$ subject to the relations
$$(T_s-v^{L(s)})(T_s+v^{-L(s)})=0 \text{ and }
\underbrace{T_sT_tT_s\cdots}_{m_{st}\ \text{terms}}=
\underbrace{T_tT_sT_t\cdots}_{m_{st}\ \text{terms}}
\ \text{ for }s,t\in S \text{ with } m_{st}{\not=}\infty$$
where $m_{st}$ is the order of $st$.
Given $w\in W$, we put
$T_w=T_{s_1}\cdots T_{s_n}$, where $w=s_1\cdots s_n$ is a reduced
decomposition.

\medskip
Let $i$ be the ring involution of $H$ given by $i(v^{\gamma})=v^{-\gamma}$ for
$\gamma\in\Gamma$
and $i(T_s)=T_s^{-1}$.
We denote by $\{C_w\}_{w\in W}$ the Kazhdan-Lusztig basis of $H$.
It is uniquely defined by the properties that $i(C_w)=C_w$ and
$C_w-T_w\in \bigoplus_{w'\in W} \BZ[\Gamma_{<0}]T_{w'}$.

\medskip
We introduce the partial order $\prec_L$ on $W$. It is the transitive closure
of the relation given by $w'\prec_L w$ if there is $s\in S$ such that
the coefficient of $C_{w'}$ in the decomposition of $C_sC_w$ in the
Kazhdan-Lusztig basis is non-zero. We define $w\sim_L w'$ to be the
corresponding equivalence relation: $w\sim_L w'$ if and only if
$w\prec_L w'$ and $w'\prec_L w$. The equivalence classes are the left cells.
We define
$\prec_{LR}$ as the partial order generated by $w\prec_{LR}w'$ if
$w\prec_L w'$ or $w^{-1}\prec_L w^{\prime -1}$. As above, we define
an associated equivalence relation
$\sim_{LR}$. Its equivalence classes are the two-sided cells.

\medskip
When $\Gamma=\BZ$,
$L=l$, and $W$ is a Weyl group, a definition of left cells based on
primitive ideals in enveloping algebras was
proposed by Joseph \cite{Jo}: let $\Gg$ be a complex semi-simple Lie algebra
with Weyl group $W$. Let $\rho$ be the half-sum of the positive roots. Given
$w\in W$, let $I_w$ be the annihilator in $U(\Gg)$ of the simple module
with highest weight $-w(\rho)-\rho$. Then, $w$ and $w'$ are in the same
left cell if and only if $I_w=I_{w'}$.

\subsection{Representations and families}
Let $\Gamma$ be a left cell. Let $W_{\le\Gamma}$ (resp. $W_{<\Gamma}$)
be the set of $w\in W$ such that there is $w'\in\Gamma$ with
$w\prec_L w'$ (resp. $w\prec_L w'$ and $w{\not\in}\Gamma$).
The left cell representation of $W$ over $\BC$
associated with $\Gamma$ \cite{KaLu,Lu4} is the unique
representation, up to isomorphism, that deforms into 
the left $H$-module 
$$\Bigl(\bigoplus_{w\in W_{\le\Gamma}}\BZ[\Gamma]C_w\Bigr)/
\Bigl(\bigoplus_{w\in W_{<\Gamma}}\BZ[\Gamma]C_w\Bigr).$$

\medskip
Lusztig \cite{Lu1,Lu4} has defined the set of constructible characters of $W$
inductively as
the smallest set of characters with the following properties:
it contains the trivial character, it is stable under 
tensoring by the sign representation and it is stable under
$J$-induction from a parabolic subgroup.
Lusztig's families are the equivalences classes of irreducible characters
of $W$ for the relation generated by $\chi\sim\chi'$ if $\chi$ and
$\chi'$ occur in the same constructible character. Lusztig has
determined constructible characters and families for all $W$ and all
parameters.

Lusztig has shown for equal parameters, and conjectured in general,
that the set of left cell characters coincides with the set of constructible
characters.

\subsection{A conjecture}
Let $c\in\BR^{\CS/\sim}$. Let $\Gamma$ be the subgroup of $\BR$ generated
by $\BZ$ and $\{c_s\}_{s\in\CS}$. We endow it with the natural order on $\BR$.
Let $L:W\to\Gamma$ be the weight
function determined by $L(s)=c_s$ if $s\in S$.

\smallskip
The following conjecture is due to Gordon and Martino \cite{GoMa}.
A similar conjecture has been proposed
independently by the second author\footnote{Talk at the Enveloping algebra
seminar, Paris, December 2004.}. It is known to hold for types
$A_n$, $B_n$, $D_n$ and $I_2(n)$ \cite{Go2,GoMa,Be,Ma1,Ma2}.

\begin{conj}
The Calogero-Moser families of irreducible characters of $W$
coincide with the Lusztig families.
\end{conj}

We propose now a conjecture involving partitions of elements of $W$, via
ramification. The part dealing with left cell characters could be stated in
a weaker way, using $Q$ and not $R$, and thus not needing the choice of
prime ideals, by involving constructible characters.

\begin{conj}
\label{co:cells}
There is a choice of $\Gr_c^{\mathrm{left}}\subset\Gr_c$ such that
\begin{itemize}
\item the
Calogero-Moser two-sided
cells (resp. left cells) coincide with the Kazhdan-Lusztig two-sided cells
(resp. left cells)
\item the representations $\sum_{E\in\Irr(W)}m_{\Gamma}(E)E$, where $\Gamma$
is a Calogero-Moser left cell, coincide
with the left cell representations of Kazhdan-Lusztig.
\end{itemize}
\end{conj}

Various particular cases and general results supporting
Conjecture \ref{co:cells} are provided in \cite{BoRou}. In particular,
the conjecture holds for $W=B_2$, for all choices of parameters.

\end{document}